\documentclass[10pt]{amsart}

\usepackage{amsfonts,amsmath,amsthm,amssymb}
\usepackage{epsfig,color}
\usepackage[mathscr]{eucal}

\newtheorem{theorem}{Theorem}

\newtheorem{proposition}{Proposition}

\theoremstyle{definition}

\newtheorem{definition}{Definition}

\newcommand{\R}{\mathbb{R}}

\newcommand{\msB}{\mathscr{B}}

\usepackage[T1]{fontenc}
\usepackage{cmbright}

\begin{document}

\title{Bernstein Polynomials and $n$-Copulas}         
 
\author{M. D. Taylor}
\address{Department of Mathematics \\
	University of Central Florida \\
	Orlando, FL 32816-1364}
\email{mtaylor@pegasus.cc.ucf.edu}

\date{\today}          
% Enter your date or\today between curly braces

\begin{abstract}
	We give derivations of some basic results for the Bernstein approximation in $n$ variables that are useful in investigating copulas.  It is shown that Bernstein approximations of copulas are again copulas.  We exhibit a stochastic interpretation for the Bernstein approximation of a copula.
\end{abstract}
 
\maketitle

\allowdisplaybreaks

\section{Introduction}

Bernstein approximations of 2-copulas were introduced and studied in \cite{k97} and \cite{Li98}.  We assume the reader is familiar with copulas; see, for example, \cite{schweizsklar83} or \cite{nelsen2006}.  

This note was written to clarify for myself and my colleagues certain properties of Bernstein approximations that are useful in investigating copulas.  We derive some of the basic properties of the Bernstein approximation for functions of $n$ variables and then show that the Bernstein approximation of a copula is again a copula.  Our most significant result is a stochastic interpretation of the Bernstein approximation of a copula.  This interpretation was communicated to us by J.\ H.\ B.\ Kemperman in \cite{kemperman} for 2-copulas and we are not aware of its publication elsewhere.  

The encouragement and contributions of our colleagues P. Mikusi\'nski and X. Li to this note were crucial.

\section{Bernstein polynomials and approximations}

It is our convention that $I = [0,1]$.

\begin{definition}
	The $m$-th degree Bernstein polynomial $b_{i,m}:I \to \R$ is given by 
\[
	b_{i,m}(t) = \binom{m}{i} \, t^i (1- t)^{m-i},
\]
$i = 0,1,\ldots,m$.  We extend this to $B^n_{i,m}:I^n \to \R$ by taking $i$ to be a multi-index, $i = (i_1,\ldots,i_n)$, where each $i_k \in \{0,1,\ldots,m\}$ and setting
\[
	B^n_{i,m}(x) = b_{i_1,m}(x_1) \, b_{i_2,m}(x_2) \cdots b_{i_n,m}(x_n)
\]
where $x = (x_1, \ldots, x_n) \in I^n$.
\end{definition}

Notice that $\{B^n_{i,m}\}$ is a partition of unity over $I^n$.

Here is the intuition behind the Bernstein polynomial:  Consider the act of tossing a coin $m$ times with probability of heads on each toss being $x$.  This scenario can be represented by a random vector $X: \Omega \to \{0,1\}^m$ with the property that if $X(\omega) = (x_1,\ldots,x_m)$, then $P(x_i = 1) = x$.  We then introduce the random variable $Y$ defined by
\[
	Y(\omega) = \sum \{x_i \, : \, x_i = 1\},
\]
in other words, the number of heads that were tossed.  This is familiarly described as a binomially distributed random variable with parameters $m$ and $x$.  It is easily shown that
\[
	E(Y) = mx \quad \text{and} \quad Var(Y) = E((Y - E(Y))^2) = mx(1-x).
\]
Notice that
\[
	b_{i,m}(x) = P(Y = i).
\]

\begin{definition}
	If $f:I^n \to \R$, we define the Bernstein approximation to $f$ to be 
\[
	\msB^n_m(f) = \sum_i f(\tfrac{i}{m}) \, B^n_{i,m}
\]
where $i$ ranges over all multi-indices $i = (i_1,\ldots,i_n)$ such that each $i_k \in \{0,1,\ldots,m\}$, and by $\frac{i}{m}$ we mean the vector $\left( \frac{i_1}{m}, \ldots, \frac{i_n}{m}\right)$.
\end{definition}

It can be shown by induction that
\[
	\int_{I^n} B^n_{i,m} \, d\lambda^n = \frac{1}{(m+1)^{n}}
\]
where $\lambda^n$ is Lebesgue measure on $\R^n$.

\section{The Weierstrass approximation theorem via Bernstein polynomials}

\begin{theorem}
	If $f:I^n \to \R$ is continuous, then $\msB^n_m(f) \to f$ uniformly on $I^n$ as $m \to \infty$.
\end{theorem}

\begin{proof}
Choose $\epsilon > 0$.  Since $f$ is uniformly continuous on $I^n$, there exists $\delta > 0$ with the property that if $x = (x_1, \ldots, x_n)$, $y = (y_1, \ldots, y_n)$, and $|x_i - y_i|<\delta$ for all $i$, then $|f(x) - f(y)| < \epsilon$.  In what follows, it is convenient to define $d$ by $d(x,y) = \max \{|x_1 - y_1|, \ldots, |x_n - y_n|\}$.

Set $f_m = \msB_m^n(f)$.  We suppose that $m$ is so large that $\frac{1}{4 m \delta^2} < \epsilon$; we shall show that this makes $|f_m - f|$ ``small.''  Choose $x = (x_1, \ldots, x_n) \in I^n$.  Then
\[
	|f_m(x) - f(x)| \leq \underset{d\left( \frac{i}{m},x\right) < \delta}{\sum} |f(\tfrac{i}{m})-f(x)| \, B^n_{i,m}(x) + \underset{d\left( \frac{i}{m},x\right) \geq \delta}{\sum} |f(\tfrac{i}{m})-f(x)| \, B^n_{i,m}(x)
\]
where $i = (i_1,\ldots,i_n)$, a multi-index.  We see that 
\[
	\underset{d\left( \frac{i}{m},x\right) < \delta}{\sum} |f(\tfrac{i}{m})-f(x)| \, B^n_{i,m}(x) < \epsilon
\]
by uniform continuity of $f$.  To find a bound for the other term, we first introduce independent, binomially distributed random variables $X_1, \ldots, X_n$ with parameters $m$ and $x$.  By Tchebycheff's inequality, for each $j = 0,1,\ldots,m$ we have 
\[
	P\left( \left| \frac{X_j}{m} - x_j \right| \geq \delta \right) \leq \frac{x_j (1-x_j)}{m \, \delta^2} \leq 
		\frac{1}{4 \, m \, \delta^2} < \epsilon.
\]
Let $M = \max f$.  Then 
\[
	 \underset{d\left( \frac{i}{m},x\right) \geq \delta}{\sum} |f(\tfrac{i}{m})-f(x)| \, B^n_{i,m}(x)\; \leq 
	 	\; 2  M  \underset{d\left( \frac{i}{m},x\right) \geq \delta}{\sum}  B^n_{i,m}(x).
\]
We see that 
\begin{gather*}
	\underset{d\left( \frac{i}{m},x\right) \geq \delta}{\sum}  B^n_{i,m}(x)  =  
		\underset{d\left( \frac{i}{m},x\right) \geq \delta}{\sum} 
		P\left(\left( \frac{X_1}{m},\ldots,\frac{X_n}{m}\right) = \frac{i}{m}\right) \\
	\leq \sum_{j=1}^n P\left( \left| \frac{X_j}{m} - x_j \right| \geq \delta \right)  < n \epsilon.
\end{gather*}
Thus $|f_m(x) - f(x)| < \epsilon + 2Mn\epsilon$, and we are done.
\end{proof}

\section{Derivatives of Bernstein approximations}

\subsection{Derivatives of Bernstein polynomials}

If $f : I^n \to \R$, set 
\begin{equation} \label{fm_form}
	f_m = \msB_m^n(f) = \sum_i f \left(\tfrac{i}{m}\right) \, B_{i,m}^n = 
		  \sum_i f \left(\tfrac{i}{m}\right) \, b_{i_1,m} \otimes \cdots \otimes b_{i_n,m}
\end{equation}
where $i = (i_1,\ldots,i_n)$ and $i_k = 0,1,\ldots,m$ and the symbolism $g \otimes h$ is interpreted to mean $(g \otimes h)(u,v) = g(u) \, h(v)$.  We want to compute partial derivatives of $f_m$.  In particular we want compute the mixed partial $\frac{\partial^n f_m}{\partial x_1 \cdots \partial x_n}$ which we denote $\partial f_m$.  In the case where $f_m$ is a cumulative probability distribution function, $\partial f_m$ is the associated probability density.

It is convenient at this point to introduce another notation.  Let $g : A \to \R$ where $A \subseteq \R^n$.  Suppose $v \in \R^n$ such that $A \cap (A-v) \ne \emptyset$.  We then define a function $\Delta_v g: A \cap (A-v) \to \R$ by $\Delta_v g(p) = g(p+v)-g(p)$.  That is, $\Delta_v g(p)$ is the variation of $g$ starting at $p$ in the direction $v$.  Next, let $e_1, \ldots, e_n$ be the standard orthonormal basis for $\R^n$, that is, $e_1 = (1,0,\ldots,0)$, $e_2 = (0,1,0,\ldots,0)$, etc.  For $f : I^n \to \R$ and $i = (i_1, \ldots, i_n)$, where $i_k = 0,1,\ldots,m-1$, we define 
\[
	\Delta^n_{i,m}f = \Delta_{\frac{1}{m} e_1} \Delta_{\frac{1}{m} e_2} \cdots 
		\Delta_{\frac{1}{m} e_n} f \left( \tfrac{i}{m} \right).
\]
We can think of $\Delta^n_{i,m}f$ as the variation of $f$ over the $n$-dimensional square $\left[ \tfrac{i_1}{m}, \tfrac{i_1+1}{m} \right) \times \cdots \times \left[ \tfrac{i_n}{m}, \tfrac{i_n+1}{m} \right)$.

Returning to the problem of derivatives, one calculates
\[
	b_{i,m}' = \begin{cases}
		-m \, b_{0,m-1}, \quad &i = 0, \\
		m \, (b_{i-1,m-1} - b_{i,m-1}), \quad &0 < i < m , \\
		m \, b_{m-1,m-1}, \quad &i = m.
	\end{cases}
\]
If we set $b_{-1,m-1} = b_{m,m-1} = 0$, then we may reduce this to 
\begin{equation} \label{bern_deriv}
	b_{i,m}' = m \, (b_{i-1,m-1} - b_{i,m-1}), \quad  0 \leq i \leq m.
\end{equation}
If one then considers the case where $n=1$ so that 
\[
	f_m = \sum_{i=0}^m f \left( \tfrac{i}{m} \right) \, b_{i,m},
\]
where $i$ is now an integer, then one easily calculates
\begin{equation} \label{fm_prime}
	f_m' = m \, \sum_{j=0}^{m-1} \left( \Delta_{\frac{1}{m} e_1}  \, 
		f \left( \tfrac{j}{m} \right) \right) \, b_{j,m-1}.
\end{equation}
We then pass to the general $n$-dimensional case where $f_m$ has the form given in (\ref{fm_form}) and by repeatedly invoking the 1-dimensional case and Equation (\ref{fm_prime}), we obtain 
\[
	\partial f_m = m^n \, \sum_j \left( \Delta^n_{j,m} f  \right) \, 
		b_{j_1,m-1} \otimes \cdots \otimes b_{j_n,m-1}
\]
where $j = (j_1, \ldots, j_n)$ and $j_k = 0,1, \ldots, m-1$.

It is well-known that the Bernstein approximation of a copula is again a copula (see, for example, \cite{k97} and \cite{Li98}), but this is also an immediate consequence of this last formula:

\begin{theorem}
	The Bernstein approximation of an $n$-copula is again an $n$-copula.
\end{theorem}

\begin{proof}
	Let $C_m = \msB^n_m(C)$ where $C$ is an $n$-copula.  The boundary conditions for a copula are easily checked.  The only questionable condition is whether or not $C_m$ is $n$-increasing.  But this follows from the fact that the terms of 
\[
	\partial C_m = m^n \, \sum_j \left( \Delta^n_{j,m} C  \right) \, b_{j_1,m-1} \otimes \cdots \otimes b_{j_n,m-1}
\]
are nonnegative.  
\end{proof}

\subsection{Some identities and estimates}

In what follows, we assume that $x \in I$, $m = 1, 2, 3, \ldots$, and $i = 0, 1, \ldots, m$.

The following is straightforward to establish by induction over $i$:

\begin{proposition} \label{prop_a}
\[
	\left( \frac{i}{m} - x \right) \, b_{i,m}(x) = x \, (1-x) \, (b_{i-1, m-1}(x) - b_{i,m-1}(x) ).
\]
\end{proposition}

\begin{proposition}
\[
	\sum_{i=0}^m \left| x - \frac{i}{m} \right| \, | b_{i-1,m-1}(x) - b_{i,m-1}(x) | = \frac{1}{m}.
\]
\end{proposition}

\begin{proof}
	We assume that $X$ is a binomially distributed random variable with parameters $x$ and $m$ and make use of Proposition \ref{prop_a}:
\begin{gather*}
	\sum_{i=0}^m \left| x - \frac{i}{n} \right| \, | b_{i-1,m-1}(x) - b_{i,m-1}(x) | = 
	  \frac{1}{x (1-x)} \, \sum_{i=0}^m \left( x - \frac{i}{m} \right)^2 \, b_{i,m}(x) \\
	= \frac{1}{x (1-x)} \, \frac{1}{m^2} \, Var (X) = \frac{1}{m}.
\end{gather*}
\end{proof}

\begin{proposition} 
\[
	\sum_{i=0}^m \left| \frac{i}{m} - x \right| \, b_{i,m}(x) = 
	2 \, x \, (1-x) \, b_{i_0,m-1}(x)
\]
where $i_0 = \lfloor mx \rfloor$, the greatest integer less than or equal to $mx$.
\end{proposition}

\begin{proof}
	Let us assume $x$ is irrational, $0 < x < 1$.  There is a unique nonnegative integer, namely $i_0 = \lfloor mx \rfloor$, such that
\[
	\frac{i_0}{m} < x < \frac{i_0 + 1}{m}.
\]
We then perform a calculation in which we invoke Proposition \ref{prop_a}:
\begin{gather*}
	\sum_{i=0}^m \left| \frac{i}{m} - x \right| \, b_{i,m}(x) = 
		\sum_{i=0}^{i_0} \left( x - \frac{i}{m} \right) \, b_{i,m}(x) + 
		\sum_{i=i_0+1}^{m} \left( \frac{i}{m} - x \right) \, b_{i,m}(x) \\
	= x \, (1-x) \, \left( \sum_{i=0}^{i_0} (b_{i,m-1}(x) - b_{i-1,m-1}(x)) + 
		\sum_{i=i_0+1}^{m} (b_{i-1,m-1}(x) - b_{i,m-1}(x)) \right) \\
	= 2 \, x \, (1-x) \, b_{i_0,m-1}(x).
\end{gather*}
We then obtain the proof for general $x$ by invoking continuity.
\end{proof}

A proof of the following has been shown to us informally by our colleague Xin Li, but it can also be found on p. 15 of \cite{lorentz_1986}.

\begin{proposition} \label{prop_b}
\[
	\sum_{i=0}^m \left| \frac{i}{m} - x \right| \, b_{i,m}(x) = O\left( \frac{1}{\sqrt{m}} \right).
\]
\end{proposition}

\begin{proposition} \label{prop_c}
	For every $\delta > 0$ and $x \in I$,
\[
	\sum_{|x - i/m| \geq \delta} b_{i,m}(x) = o \left( \frac{1}{m} \right).
\]
\end{proposition}

\begin{proof}
	From \cite[p. 304]{devore_1993}, we find that for each $\delta > 0$ and $s = 1, 2, \ldots$, there exists $C = C(\delta , s)$ such that
\[
	\sum_{|x - i/m| \geq \delta} b_{i,m}(x) \leq C \, m^{-s}
\]
for $m = 1, 2, \ldots$ and $x \in I$.
\end{proof}

\subsection{Convergence of first derivatives of Bernstein approximations}

\begin{theorem}
	Suppose that $f : I^n \to \R$ is a bounded function.  Then for all $x = (x_1, \ldots, x_n) \in (0,1)^n$ at which $f$ is differentiable and for $k = 1, \ldots, n$, we have 
\[
	\underset{m \to \infty}{\lim} \frac{\partial \msB_m^n(f)}{\partial x_k} (x) = 
		\frac{\partial f}{\partial x_k} (x).		
\]
\end{theorem}

\begin{proof}
	Let us set 
\[
	f_m = \msB_m^n(f) = \sum_i f \left(\tfrac{i}{m} \right) \, b_{i_1,m} \otimes \cdots \otimes b_{i_n,m}
\]
where $i = (i_1, \ldots, i_n)$, each $i_k \in \{0, 1, \ldots, m\}$, and
\[
	\sum_i = \sum_{i_1=0}^m \sum_{i_2=0}^m \cdots \sum_{i_n=0}^m.
\]  
We prove the proposition for the case $k=1$.

First, we have
\begin{equation} \label{bern_deriv}
	\frac{\partial f_m}{\partial x_1} = m \, \sum_i f(\tfrac{i}{m}) \, (b_{i_1-1,m-1} - b_{i_1,m-1}) \otimes b_{i_2,m} \otimes \cdots \otimes b_{i_n,m}.
\end{equation}
Second, by the differentiability of $f$ at $x$, we have
\begin{equation} \label{f_deriv}
	f(\tfrac{i}{m}) = f(x) + \sum_{k=1}^n \frac{\partial f}{\partial x_k} (x) \, 
	\left( \frac{i_k}{m} - x_k \right) + \eta( \tfrac{i}{m} - x ) \, \left| \frac{i}{m} - x \right|
\end{equation}
where 
\[
	\left| \frac{i}{m} - x \right| = \sqrt{ \sum_{k=1}^n \left( \frac{i_k}{m} - x_k \right)^2}
\]
and $\eta(s) \to 0$ as $s \to 0$ in $\R^n$.  Next, making use of (\ref{f_deriv}), it can be shown there is a constant $M$, dependent on $x$ and $n$ but independent of $m$, such that $| \eta(s) | \leq M$.  This can be done by considering the case where $s$ is ``close'' to 0 and the case where $s$ is some fixed distance from 0.

Next, substituting from (\ref{f_deriv}) into (\ref{bern_deriv}), we see that $(\partial f_m / \partial x_1)(x)$ becomes 
\begin{align*}
	m \, \sum_i \Bigg(   f(x) +& \sum_{k=1}^n \frac{\partial f}{\partial x_k} (x) \, 	\left( \frac{i_k}{m} - x_k \right) + \eta( \tfrac{i}{m} - x ) \, \left| \frac{i}{m} - x \right| \Bigg) \\
&(b_{i_1-1,m-1}(x_1) - b_{i_1,m-1}(x_1)) \, b_{i_2,m}(x_2) \, \cdots  b_{i_n,m}(x_n) 
\end{align*}
Now $\sum_{i_1} (b_{i_1-1,m-1}(x_1) - b_{i_1,m-1}(x_1)) = b_{-1,m-1}(x_1) - b_{m,m-1}(x_1) = 0$.  Thus 
\[
	\sum_i  f(x) \big(b_{i_1-1,m-1}(x_1) - b_{i_1,m-1}(x_1) \big) \, b_{i_2,m}(x_2) \, \cdots  b_{i_n,m}(x_n) = 0,
\]
and for every $k>1$ we have
\begin{gather*}
	\sum_i \frac{\partial f}{\partial x_k} (x) \, \left( \frac{i_k}{m} - x_k \right) 
	\big(b_{i_1-1,m-1}(x_1) - b_{i_1,m-1}(x_1) \big) \, b_{i_2,m}(x_2) \, \cdots  b_{i_n,m}(x_n) \\ 
	= 0.
\end{gather*}
On the other hand, it is easily seen that
\begin{align*}
	m \sum_{i_1=0}^m & \left( \frac{i_1}{m} - x_1 \right) \, (b_{i_1-1,m-1}(x_1) - b_{i_1,m-1}(x_1)) \\ 
	& = \sum_{i_1} i_1 \, (b_{i_1-1,m-1}(x_1) - b_{i_1,m-1}(x_1)) = 1,
\end{align*}
and we know that $\sum_{i_k} b_{i_k,m}(x_k) = 1$ for $k>1$, therefore
\begin{gather*}
	m \sum_i \frac{\partial f}{\partial x_1}(x) \, \left( \frac{i_1}{m} - x_1 \right)  \, \big(b_{i_1-1,m-1}(x_1) - b_{i_1,m-1}(x_1) \big) \, b_{i_2,m}(x_2) \, \cdots  b_{i_n,m}(x_n) \\
	= \frac{\partial f}{\partial x_1}(x).
\end{gather*}
Thus we can write
\[
	\frac{\partial f_m}{\partial x_1}(x) = \frac{\partial f}{\partial x_1}(x) + S_m	
\]
where
\[
	S_m = m \sum_i \eta( \tfrac{i}{m} - x ) \, \left| \frac{i}{m} - x \right| \, \big(b_{i_1-1,m-1}(x_1) - b_{i_1,m-1}(x_1) \big) \, b_{i_2,m}(x_2) \, \cdots  b_{i_n,m}(x_n).
\]

Our task now reduces to showing $S_m \to 0$.  Choose $\epsilon > 0$.  There exists $\delta > 0$ such that if $|s|<\delta$, then $|\eta(s)|<\epsilon$.  Let us set $\delta_i = |(i/m)-x|$ and then break $S_m$ into two pieces, $S_m = S_m^> + S_m^{\geq}$, where
\[
	S_m^< = \sum_{\delta_i < \delta} \quad \text{and} \quad S_m^{\geq} = \sum_{\delta_i \geq \delta}.
\]

We first consider $S_m^<$.  By Proposition \ref{prop_a}, 
\[
	b_{i_1-1,m-1}(x_1) - b_{i_1,m-1}(x_1) = \frac{ \frac{i_1}{m} - x_1}{x_1(1-x_1)} \, b_{i_1,m}(x_1).
\]
Then
\begin{gather*}
	| S_m^<| \leq m \, \epsilon \sum_{\delta_i < \delta} \delta_i \, \big(b_{i_1-1,m-1}(x_1) - b_{i_1,m-1}(x_1) \big) \, b_{i_2,m}(x_2) \, \cdots  b_{i_n,m}(x_n) \\
	\leq \frac{m \, \epsilon}{x_1(1-x_1)} \sum_{\delta_i < \delta} \delta_i \, \left| \frac{i_1}{m}-x_1 \right| \, b_{i_1,m}(x_1) \, b_{i_2,m}(x_2) \cdots b_{i_n,m}(x_n) \\
	\leq \frac{m \, \epsilon}{x_1(1-x_1)} \sum_{i} \sum_{k=1}^n \left| \frac{i_k}{m}-x_k \right| \, \left| \frac{i_1}{m}-x_1 \right| \, b_{i_1,m}(x_1) \cdots b_{i_n,m}(x_n).
\end{gather*}
Now $\sum_{i_1}(i_1-mx_1)^2 b_{i_1,m}(x_1)$ is the variance of a binomially distributed random variable, so
\[
	\sum_{i_1} \left( \frac{i_1}{m} - x_1 \right)^2 b_{i_1,m}(x_1) = \frac{x_1(1-x_1)}{m}.
\]
On the other hand, for $k \ne 1$, we have by Proposition \ref{prop_b},
\[
	\sum_{i_1} \sum_{i_k} \left| \frac{i_1}{m} - x_1 \right| \, \left| \frac{i_k}{m} - x_k \right| \, b_{i_1,m}(x_1) \, b_{i_k,m}(x_k) = O\left( \tfrac{1}{\sqrt{m}} \right) \, O\left( \tfrac{1}{\sqrt{m}} \right) = O( \tfrac{1}{m} ).
\]
It follows that 
\[
	|S_m^<| \leq \frac{ m \, \epsilon}{x_1(1-x_1)} \, \frac{x_1(1-x_1)}{m} + 
	m \, (n-1) \, \epsilon \, O( \tfrac{1}{m} ) = \epsilon \, O(1).
\]

Now we turn to $S_m^{\geq}$.  We know the following:
\begin{enumerate}

	\item  $| \eta | \leq M$.

	\item  $\delta_i \leq \sqrt{n}$.

	\item  $\left\{ i \, : \, \left| \frac{i}{m} - x \right| \geq \delta \right\} \subseteq \bigcup_{k=1}^n \left\{ i \, : \, \left| \frac{i_k}{m} - x_k \right| \geq \frac{\delta}{n} \right\}$.
\end{enumerate}
Using these facts plus Proposition \ref{prop_c}, we obtain
\begin{gather*}
	|S_m^{\geq}| \leq m \, M \, \sqrt{n} \sum_{ \left| \frac{i}{m}-x \right| \geq \delta } b_{i_1,m}(x_1) \cdots b_{i_n,m}(x_n) \\
	\leq m \, M \, \sqrt{n} \sum_{ \left| \frac{i_k}{m}-x_k \right| \geq \frac{\delta}{n} } \; \sum_{k=1}^n b_{i_k,m}(x_k) \\
	= m \, M \, \sqrt{n} \; \; n \; o ( \tfrac{1}{m} ) = o(1).
\end{gather*}

We therefore conclude that 
\[
	|S_m| \leq \epsilon \, O(1) + o(1) \to 0
\]
as $m \to \infty$.
\end{proof}

\section{A probabilistic interpretation of the Bernstein approximation of a copula}

It is our goal here to construct random variables such that the Bernstein approximation is the cumulative distribution function of these new random variables.  This probabilistic interpretation was brought to our attention by J. H. B. Kemperman in \cite{kemperman}.  

Let $C$ be an $n$-copula.  Suppose it is the cumulative distribution function of the ordered $n$-tuple of random variables $(X_1,\ldots,X_n)$ where each $X_i$ is uniformly distributed over $I$.  Let $m$ be a ``large'' natural number and $C_m$ be the $m \times \cdots \times m$ Bernstein approximation of $C$; that is 
\[
	C_m(x_1, \ldots, x_n) = \sum_{i_1, \ldots, i_n = 0}^m C\left( \tfrac{i_1}{m}, \ldots, \tfrac{i_n}{m} \right) \; b_{i_1,m}(x_1) \cdots b_{i_n,m}(x_n)
\]
where $x_1, \ldots, x_n \in I$.  

Next, for $i = 1, \ldots, n$ and $j = 1, \ldots, m$, we let $X_i^j$ be independent random variables that are uniformly distributed over $I$ and have the property that $(X_1, \ldots, X_n)$ and $X_i^j$ are independent for all $i,j$.  If it is helpful, we may regard these random variables as being defined over the space $I^n \times I^{mn}$ and having probability measure $P$ where $P$ has the form $\mu_C \times \lambda^{mn}$ and where it is understood that $\mu_C$ is the probability measure induced on $I^n$ by $C$ and $\lambda^{mn}$ is Lebesgue measure on $I^{mn}$.   

Now for each $i$, let $X_i^{(1)}, \ldots, X_i^{(m)}$ be the order statistics for $X_i^{1}, \ldots, X_i^{m}$.  That is, whenever 
\[
	X_i^{j_1} < \cdots < X_i^{j_m},
\]
then $X_i^{(k)} = X_i^{j_k}$.  

It may be helpful to notice that by the independence and uniform distribution of our original random variables, for all $i,j,r,s$ and all $x \in I$ we have 
\begin{gather*}
	P(X_i^j = X_r^s) = 0 \quad \text{if} \quad X_i^j \ne X_r^s, \\
	P(X_i^{(j)} = x) = 0, \\
	P(X_i^{(j)} = X_r^{(s)}) = 0 \quad \text{if} \quad X_i^{(j)} \ne X_r^{(s)}.
\end{gather*}

Next, if $k$ is the multi-index $(k_1, \ldots, k_n)$ where $k_r = 0, 1, \ldots, m-1$, then we define 
\[
	I^n_{m,k} = 
	\left( \frac{k_1}{m}, \frac{k_1+1}{m} \right) \times \cdots \times \left( \frac{k_n}{m}, \frac{k_n+1}{m}  \right) = 
	\left\{ \frac{k_r}{m} < X_r < \frac{k_r+1}{m} \, : \, r = 1, \ldots, n \right\}.
\]
This is a slight abuse of notation since $X_r$ ``lives'' in $I^n \times I^{mn}$ rather than $I^n$, however the reader should easily make whatever mental adjustments are necessary in the arguments that follow.  We then take $\chi_k$ to be the characteristic function of $I^n_{m,k}$ with the understanding that the domain of $\chi_k$ is $I^n$.

Finally we define 
\[
	Y_r = \sum_{k_1, \ldots, k_n=0}^{m-1} \chi_k(X_1, \ldots, X_n) \; X_r^{(k_r)}
\]
where $r=1, \ldots, n$ and $k = (k_1, \ldots, k_n)$.  We then have the following:

\begin{theorem}
	$C_m(x_1, \ldots, x_n) = P( Y_1 < x_1, \ldots, Y_n < x_n )$ where $(x_1, \ldots, x_n) \in I^n$.
\end{theorem}

\begin{proof}
	For $x \in I$ and $i = 1, \ldots, n$, we form a random variable $X_i^*(x)$ by setting 
\[
	X_i^*(x) = \text{ the number of $X_i^1, \ldots, X_i^m$ less than $x$.}
\]
We see that the following are true:
\begin{enumerate}
	\item  $X_i^*(x)$ takes values in $\{0,1, \ldots, m\}$.
	\item  $X_i^*(x)$ and $X_j^*(y)$ are independent for $i \ne j$.
	\item  $P(X_i^*(x) = k) = \binom{m}{k} \, x^k \, (1-x)^{m-k} = b_{i,m}(x)$.
\end{enumerate}
We also see that 
\begin{equation} \label{ordstatcount}
	\begin{split}
	\{ X_i^*(x)=0 \} &\overset{P\text{-a.e.}}{=} \{x < X_i^{(1)} \}, \\
	\{ X_i^*(x)=k \} &\overset{P\text{-a.e.}}{=} \{ X_i^{(k)} < x < X_i^{(k+1)} \} \quad \text{for } k=1, \ldots, m-1, \\
	\{ X_i^*(x)=m \} &\overset{P\text{-a.e.}}{=} \{ X_i^{(m)} < x \}.
	\end{split}
\end{equation}
From (\ref{ordstatcount}) we can easily deduce 
\begin{equation} \label{ordstatcount2}
	\{ X_i^*(x) \geq k \} \overset{P\text{-a.e.}}{=} \{ X_i^{(k)} < x \}
\end{equation}
for $i = 1, \ldots, n$ and $k = 0,1, \ldots, m$.

We now consider the Bernstein approximation to $C$.  Let $x_1, \ldots, x_n \in I$.  Then 
\begin{align*}
	C_m(x_1, &\ldots, x_n) = \sum_{i_1, \ldots, i_n = 1}^m C\left( \tfrac{i_1}{m}, \ldots, \tfrac{i_n}{m}\right) \, b_{i_1,m}(x_1) \ldots b_{i_n,m}(x_n) \\
	=& \sum_{i_1, \ldots, i_n = 1}^m \bigg[ \sum_{k_1=1}^{i_1} \cdots \sum_{k_n=1}^{i_n} P \left( \frac{k_r-1}{m} < X_r < \frac{k_r}{m} \; : \; r = 1, \ldots, n \right) \\
	& \hspace{1in} P\left( X_r^*(x_r) = i_r \, : \, r = 1, \ldots, n \right) \, \bigg].
\end{align*}
Since
\[
	\sum_{i_r=1}^m \, \sum_{k_r=1}^{i_r} = \sum_{k_r=1}^m \, \sum_{i_r=k_r}^m,
\]
we have
\begin{align*}
	C_m(x_1, &\ldots, x_n) = \\
	=& \sum_{k_1, \ldots, k_n=1}^m \; \sum_{i_1=k_1}^m \cdots \sum_{i_n=k_n}^m \bigg[ P \left( \frac{k_r-1}{m} < X_r < \frac{k_r}{m} \; : \; r = 1, \ldots, n \right) \\
	& \hspace{1in} P\left( X_r^*(x_r) = i_r \, : \, r = 1, \ldots, n \right) \, \bigg] \\
	=& \sum_{k_1, \ldots, k_n=1}^m \bigg[ P \left( \frac{k_r-1}{m} < X_r < \frac{k_r}{m} \; : \; r = 1, \ldots, n \right) \\
	& \hspace{1in} P\left( X_r^*(x_r) \geq k_r \, : \, r = 1, \ldots, n \right) \, \bigg] \\
	=& \sum_{k_1, \ldots, k_n=1}^m \bigg[ P \left( \frac{k_r-1}{m} < X_r < \frac{k_r}{m} \; : \; r = 1, \ldots, n \right) \\
	& \hspace{1in} P\left( X_r^{(k_r)} < x_r \, : \, r = 1, \ldots, n \right) \, \bigg] 
\end{align*}
where the last step follows from (\ref{ordstatcount2}).  

We now consider the cumulative distribution function of $(Y_1, \ldots, Y_n)$.
\begin{align*}
	P(Y_r &< x_r \, : \, r = 1, \ldots, n) \\
	&= \sum_{k_1, \ldots, k_n=1}^m P\left(Y_r < x_r, \; (X_1, \ldots, X_n) \in I^n_{m,k} \; : \; r = 1, \ldots, n \right) \\
	&= \sum_{k_1, \ldots, k_n=1}^m \bigg[ P\left( X_r^{(k_r)} < x_r \, : \, r = 1, \ldots, n \right) \\
	& \hspace{1in} P \left( \frac{k_r-1}{m} < X_r < \frac{k_r}{m} \; : \; r = 1, \ldots, n \right) \bigg]
\end{align*}
by the definition of $Y_r$ and the independence of the random variables.  We see from this that 
\[
	C_m(x_1, \ldots, x_n) = P(Y_1 < x_1, \ldots, Y_n < x_n).   \qedhere
\]
\end{proof}

\end{document}